\documentclass{amsart}
\usepackage{latexsym,amssymb,amsmath,amsthm,latexsym}
\usepackage{amsfonts}
\theoremstyle{definition}

\newtheorem{example}{Example}
\newtheorem{corollary}{Corollary}

\newtheorem{proposition}{Proposition}

\newtheorem{remark}{Remark}

\newtheorem{identity}{Identity}
\title{Cassini, d'Ocagne, and Vajda identities for $n$-step Fibonacci  numbers}
\author{Milan Janji\'c}
\date{\today}
\begin{document}
\maketitle
\begin{center}Department for Mathematics and Informatics, University of Banja
Luka\\Republic of Srpska, Bosnia and Herzegovina\end{center}

\begin{abstract}Results of this paper concern $n$-determinants which we defined in
the paper \cite{jan}. In the paper \cite{jabo}, $2$-determinants are  considered. In this paper, we extend results from \cite{jabo} on $n$-determinants by proving that   Cassini,
d'Ocagne, Catalan and Vajda identities may be extended to hold for  $n$-step Fibonacci numbers.
 \end{abstract}
 \section{Introduction} We restate a particular case of a result proved in
 \cite{jan}.
Let $n$ and $r$ be  positive integers. We consider the following
$n+r-1$ by $r$ matrix $P$:
\begin{displaymath}P=\begin{pmatrix}
1&0&\cdots&0&0\\
1&1&\cdots&0&0\\
\vdots&\vdots&\cdots&\vdots&\vdots\\
1&1&\cdots&1&0\\
-1&1&\cdots&1&1\\
0&-1&\cdots&1&1\\
\vdots&\vdots&\cdots&\vdots&\vdots\\
0&0&\cdots&1&1\\
0&0&\cdots&-1&1
\end{pmatrix}.\end{displaymath}
We note that in each  column we have exactly $n$ ones.
We connect the matrix $P$ with  recursively given sequence of
vector-columns in the following way: Let $A=(A_1|A_2|\ldots|A_n)$ be a square matrix
of order $n$, where $A_1,\ldots,A_n$ are columns of $A.$   We define a block matrix
$A_{n+r}$ of $n$ rows and $n+r$ columns in the following way: First $n$ columns of
$A_{n+r}$ are  $A_i,(1,2,\ldots,n)$, and
 \begin{equation}\label{r1}
 A_{n+j}=\sum_{i=j}^{n+j-1}A_i,\;(j=1,2,\ldots,r).\end{equation}

For a  sequence  $1\leq j_1<j_2<\cdots<j_r<n+r$ of positive
integers we  let $M=M(\widehat{j_1},\widehat{j_2},\ldots,\widehat{j_r})$ denote  the
minor of  $A_{n+r}$ of order $n$,  obtained  by deleting  columns
$j_1,j_2,\ldots,j_r$ of $A_{n+r}.$
Note that the last column of $A_{n+r}$ can not be deleted.
The sign of $M$  is defined as
      ${\rm sgn}(M)=(-1)^{nr+j_1+j_2+\cdots+j_r+\frac{(r-1)r}{2}}.$
We let  $Q_r=Q_r(j_1,\ldots,j_r)$ denote the sub-matrix of order $r$,
laying in rows $j_1,j_2,\ldots,j_r$ of $P.$

 Denote by $\{i_1,i_2,\ldots,i_n\}$ the
numbers of columns of $B_{n+r}$ in which the minor $M$ lies. We have $i_n=r+n$, since the last column of $A_{n,n+r}$ can not be deleted. We have
$\{i_1,i_2,\ldots,i_{n-1}\}\cup\{j_1,j_2,\ldots,j_r\}=\{1,2,\ldots,n+r-1\}.$
It follows that $M=M(\widehat{j_1},\widehat{j_2},\ldots,\widehat{j_r})=M(i_1,i_2,\ldots,i_{n-1},n+r).$
and ${\rm sgn}(M)=(-1)^{\frac{n(n-1)}{2}+i_1+i_2+\ldots,i_{n-1}}$.
In \cite{jan} the following result is proved:
\begin{proposition}\label{th1} Let $1\leq j_1<\cdots<j_r<r+n$ be a sequence of positive
integers. Then,
\begin{equation}\label{ttt}M({i_1},\ldots,i_{n-1},r+n)={\rm sgn}(M)\cdot
\det Q(j_1,j_2,\ldots,j_r)\cdot \det A.\end{equation}
\end{proposition}
We stress a particular case of this result, when $i_k=k,(k=1,2,\ldots,n-1)$.
In this case, we have $j_1=n,j_2=n+1,\ldots,j_r=n+r-1$, so that the matrix $Q_r$
has the following form
\begin{displaymath}
Q_r(n,\ldots,n+r-1)=\begin{vmatrix}1&1&\cdots&1&0&0&\cdots&0\\
-1&1&\cdots&1&1&0&\cdots&0\\
0&-1&\cdots&1&1&1&\cdots&0\\
\vdots&\vdots&\vdots&\cdots&\vdots&\vdots&\vdots&\vdots\\
0&0&0&\cdots&0&-1&1&1\\
0&0&0&\cdots&0&0&-1&1
\end{vmatrix}.
\end{displaymath}
For $r=1,2,\ldots$, we define $\det Q_r=F^{(n)}_{r}$.
It follows that
\begin{displaymath}
F^{(n)}_1=1, F^{(n)}_2=2,\ldots,F^{(n)}_{n}=2^{n-1},
\end{displaymath}
and $r>n$ we have
\begin{displaymath}
F^{(n)}_r=F^{(n)}_{r-1}+F^{(n)}_{r-2}+\cdots+F^{(n)}_{r-n}.
\end{displaymath}
Hence, the sequence $\det Q_1,\det Q_2,\ldots$ consists of Fibonacci $n$-step numbers.
It is obvious that ${\rm sgn}(M)=1$. Using equation (\ref{ttt}), we obtain
\begin{proposition}\label{pp1}[Generalized d'Ocagne identity] The following equation holds
\begin{displaymath}M(1,\ldots,n-1,r+n-1)=F^{(n)}_{r}\cdot\det A.\end{displaymath}
\end{proposition}
We see that the expression $\frac{M(1,\ldots,n-1,r+n-1)}{\det A}$ does not depend on the matrix $A$ with $\det A\not=0$. In other words, the following result holds:
\begin{proposition}Assume that $A$ and $B$ are regular matrices of order $n$, and let $A_{n+r}$ and $B_{n+r}$ be formed according the same role described above.
Let $M$ and $N$ are corresponding minors. Then
\begin{displaymath}
\frac{M(1,\ldots,n-1,r+n-1)}{\det A}=\frac{N(1,\ldots,n-1,r+n-1)}{\det B}.
\end{displaymath}
\end{proposition}
Our next goal is to extend Cassini identity of Fibonacci numbers on $n$-step
Fibonacci numbers.

\section{Cassini Identity for $n$-step Fibonacci  numbers}

We note that cases $n=2$ and $n=3$ are already considered in \cite{jan}.

For the matrix $A$, we take the following  upper triangular matrix:
\begin{displaymath}
A=\begin{pmatrix}F^{(n)}_1&F^{(n)}_2&F^{(n)}_3&\cdots&F^{(n)}_{n-1}&F^{(n)}_n\\
0&F^{(n)}_1&F^{(n)}_2&\cdots&F^{(n)}_{n-2}&F^{(n)}_{n-1}\\
\vdots&\vdots&\vdots&\cdots&\vdots\\
0&0&0&\cdots&0&F^{(n)}_{1}
\end{pmatrix}.
\end{displaymath}
Hence, $\det A=1$. We next take $i_k=r+k,(i=1,2,\ldots,n-1)$. Since
$Q_r(1,2,\ldots,r)$ is a lower triangular matrix having ones on the main diagonal, we
have $\det Q_r=1$. Also, $M=M(r+1,r+2,\ldots,r+n)$, and ${\rm sgn}M=(-1)^{(n-1)r}$.
We thus obtain
\begin{proposition}[Cassini identity for $n$-step Fibonacci numbers]
\begin{equation}\label{cas}
\begin{vmatrix}F^{(n)}_{r+1}&F^{(n)}_{r+2}&\cdots&F^{(n)}_{r+n}\\
F^{(n)}_{r}&F^{(n)}_{r+1}&\cdots&F^{(n)}_{r+n-1}\\
\vdots&\vdots&\cdots&\vdots\\
F^{(n)}_{r-n+2}&F^{(n)}_{r-n+3}&\cdots&F^{(n)}_{r+1}
\end{vmatrix}=(-1)^{(n-1)r}.
\end{equation}
\begin{remark}This result is known. It is obtained in \cite{sim} using induction.
Also a number of particular identities is given in this paper.
\end{remark}

\end{proposition}
\begin{corollary} In the case $n=2$, equation (\ref{cas}) becomes Cassini identity for
Fibonacci numbers
\begin{displaymath}
\begin{vmatrix}F_{r+1}&F_{r+2}\\
F_{r}&F_{r+1}
\end{vmatrix}=(-1)^{r}.
\end{displaymath}
\end{corollary}
\begin{corollary} In the case $n=3$, we obtain Cassini identity for Tribonacci
numbers $T_1=1,T_2=1,T_3=2,\ldots$:
\begin{displaymath}
\begin{vmatrix}T_{r+1}&T_{r+2}&T_{r+3}\\
T_{r}&T_{r+1}&T_{r+2}\\
T_{r-1}&T_r&T_{r+1}
\end{vmatrix}=1.
\end{displaymath}
\end{corollary}
\section{d'Ocagne identity for $n$-step Fibonacci numbers}
We next take that $A=
\begin{pmatrix}F^{(n)}_{r+1}&F^{(n)}_{r+2}&\cdots&F^{(n)}_{r+n}\\
F^{(n)}_{r}&F^{(n)}_{r+1}&\cdots&F^{(n)}_{r+n-1}\\
\vdots&\vdots&\cdots&\vdots\\
F^{(n)}_{r-n+2}&F^{(n)}_{r-n+3}&\cdots&F^{(n)}_{r+1}
\end{pmatrix}$. For $s>0$, we form the matrix $A_{n+s}$, such that columns $n+i,(i=n+1,\ldots,s)$ are
 sums of preceding $n$ columns. We calculate
$M=M(1,2,\ldots,s-1,s+n)$. In this case, we have
${\rm sgn}(M)=(-1)^{\frac{n(n-1)}{2}+1+2+\ldots+n-1}=1$. Using Proposition \ref{pp1}, we obtain
\begin{proposition}[D'Ocagne identity for $n$-step Fibonacci numbers]
The following formula holds
 \begin{equation}\label{okag}
\begin{vmatrix}F^{(n)}_{r+1}&F^{(n)}_{r+2}&\cdots&F^{(n)}_{r+n-1}&F^{(n)}_{r+n+s-1}\\
F^{(n)}_{r}&F^{(n)}_{r+1}&\cdots&F^{(n)}_{r+n-2}&F^{(n)}_{r+n+s-2}\\
\vdots&\vdots&\cdots&\vdots&\vdots\\
F^{(n)}_{r-n+2}&F^{(n)}_{r-n+3}&\cdots&F^{(n)}_{r}&F^{(n)}_{r+s}
\end{vmatrix}=(-1)^{(n-1)r}\cdot F^{(n)}_{s}.
 \end{equation}
\end{proposition}
\begin{example} For $n=2$, (\ref{okag}) becomes d'Ocagne identity for Fibonacci
numbers:
\begin{displaymath}
\begin{vmatrix}F_{r+1}&F_{r+s+1}\\
F_{r}&F_{r+s}\end{vmatrix}
=(-1)^r\cdot F_{s}.
\end{displaymath}
\end{example}
\begin{example} For $n=3$, from  (\ref{okag}) we obtain d'Ocagne identity for
Tribonacci numbers:
\begin{displaymath}
\begin{vmatrix}T_{r+1}&T_{r+2}&T_{r+s+2}\\
T_{r}&T_{r+1}&T_{r+s+1}\\
T_{r-1}&T_r&T_{r+s}
\end{vmatrix}=T_{s}.
\end{displaymath}
\end{example}
\section{Vajda identity for $n$-step Fibonacci numbers}
 In the determinant on the right hand side of the equation (\ref{okag}) sums of
 elements in rows are no more $n$-step Fibonacci numbers. But sums of element in
 columns are. So we consider the following matrix
  \begin{equation}\label{trokg}B=\begin{pmatrix}F^{(n)}_{r+1}&F^{(n)}_{r}&\cdots&
  F^{(n)}_{r-n+3}&F^{(n)}_{r-n+2}\\
F^{(n)}_{r+2}&F^{(n)}_{r+1}&\cdots&F^{(n)}_{r-n+4}&F^{(n)}_{r-n+3}\\
\vdots&\vdots&\cdots&\vdots&\vdots\\
F^{(n)}_{r+n-1}&F^{(n)}_{r+n-2}&\cdots&F^{(n)}_{r-1}&F^{(n)}_{r}\\
F^{(n)}_{r+n+s-1}&F^{(n)}_{r+n+s-2}&\cdots&F^{(n)}_{r+s-1}&F^{(n)}_{r+s}
\end{pmatrix},
 \end{equation}
which is simple the transpose of $A$. Hence, $\det B=(-1)^{(n-1)r}\cdot
F^{(n)}_{s}$.
Permuting columns  $1$ and $n$, $2$ and $n-1$ and so on, the matrix $B$ becomes
\begin{equation}\label{trokg}C=\begin{pmatrix}F^{(n)}_{r-n+2}&F^{(n)}_{r-n+3}&\cdots&
  F^{(n)}_{r}&F^{(n)}_{r+1}\\
F^{(n)}_{r-n+3}&F^{(n)}_{r-n+4}&\cdots&F^{(n)}_{r+1}&F^{(n)}_{r+2}\\
\vdots&\vdots&\cdots&\vdots&\vdots\\
F^{(n)}_{r}&F^{(n)}_{r+1}&\cdots&F^{(n)}_{r+n-2}&F^{(n)}_{r+n-1}\\
F^{(n)}_{r+s}&F^{(n)}_{r+s+1}&\cdots&F^{(n)}_{r+n+s-2}&F^{(n)}_{r+n+s-1}
\end{pmatrix}.
 \end{equation}
For this we need $\lfloor\frac{n}{2}\rfloor$ permutation. It follows that
\begin{displaymath}\det C=(-1)^{\lfloor\frac{n}{2}\rfloor}\cdot\det B.
\end{displaymath}

Applying (\ref{okag}), we obtain
\begin{proposition}[Vajda identity for $n$-step Fibonacci numbers]
The following formula holds
\begin{displaymath}
\begin{vmatrix}F^{(n)}_{r-n+2}&F^{(n)}_{r-n+3}&\cdots&
  F^{(n)}_{r}&F^{(n)}_{p+r}\\
F^{(n)}_{r-n+3}&F^{(n)}_{r-n+4}&\cdots&F^{(n)}_{r+1}&F^{(n)}_{p+r+1}\\
\vdots&\vdots&\cdots&\vdots&\vdots\\
F^{(n)}_{r}&F^{(n)}_{r+1}&\cdots&F^{(n)}_{r+n-2}&F^{(n)}_{p+r+n-2}\\
F^{(n)}_{q+r}&F^{(n)}_{q+r+1}&\cdots&F^{(n)}_{q+r+n-2}&F^{(n)}_{p+q+r+n-2}
\end{vmatrix}=(-1)^{(n-1)r+\lfloor\frac{n}{2}\rfloor
}\cdot F^{(n)}_{p}\cdot F^{(n)}_{q}\end{displaymath}
\end{proposition}
As a particular case, for $p=q$, we obtain
\begin{corollary}[Catalan identity for $n$-step Fibonacci numbers]
The following formula holds
\begin{displaymath}
\begin{vmatrix}F^{(n)}_{r-n+2}&F^{(n)}_{r-n+3}&\cdots&
  F^{(n)}_{r}&F^{(n)}_{p+r}\\
F^{(n)}_{r-n+3}&F^{(n)}_{r-n+4}&\cdots&F^{(n)}_{r+1}&F^{(n)}_{p+r+1}\\
\vdots&\vdots&\cdots&\vdots&\vdots\\
F^{(n)}_{r}&F^{(n)}_{r+1}&\cdots&F^{(n)}_{r+n-2}&F^{(n)}_{p+r+n-2}\\
F^{(n)}_{r+p}&F^{(n)}_{r+p+1}&\cdots&F^{(n)}_{r+p+s-2}&F^{(n)}_{2p+r+n-2}
\end{vmatrix}=(-1)^{(n-1)r+\lfloor\frac{n}{2}\rfloor
}\cdot [F^{(n)}_{p}]^2.\end{displaymath}
\end{corollary}
\begin{example}In the case $n=2$, we obtain the standard Vajda and Catalan identities for Fibonacci numbers.
We have
\begin{identity}[Vajda  identity]The following formula holds
\begin{displaymath}
\begin{vmatrix}F_{r}&F_{p+r}\\F_{r+q}&F_{p+r+q}
\end{vmatrix}=(-1)^{r+1}\cdot F_p\cdot F_q.
\end{displaymath}
\end{identity}
We also have
\begin{identity}[Catalan identity] The following formula holds
\begin{displaymath}
\begin{vmatrix}F_{r-p}&F_{r}\\F_{r}&F_{p+r}
\end{vmatrix}=(-1)^{r-p}\cdot F_p^2.
\end{displaymath}
\end{identity}

\end{example}
\begin{example}[Tribonacci numbers] In the case $n=3$, we obtain the following
identities for Tribonacci numbers $T_1=1,T_2=2,T_3=4,T_5=7,\ldots$
\begin{identity}[Vajda identity for Tribonacci  numbers] The following formula holds
\begin{displaymath}
\begin{vmatrix}T_{r-1}&T_{r}&T_{p+r}\\
T_{r}&T_{r+1}&T_{p+r+1}\\
T_{r+q}&T_{r+q+1}&T_{p+r+q+1}
\end{vmatrix}= -T_p\cdot T_q.
\end{displaymath}
\end{identity}
For $p=q$ we have
\begin{identity}[Catalan  identity for Tribonacci numbers]  We have
\begin{displaymath}
\begin{vmatrix}T_{r-1}&T_{r}&T_{p+r}\\
T_{r}&T_{r+1}&T_{p+r+1}\\
T_{r+p}&T_{r+p+1}&T_{2p+r+1}
\end{vmatrix}= -T_p^2.
\end{displaymath}
\end{identity}
\end{example}


\begin{thebibliography}{9}
\bibitem{jabo}D. Bogdani\'c, M. Janji\'c, Arithmetic of some sequences via
    $2$-determinants, Bull. Int. Math. Virtual Inst.,11(3)(2021), 517-526.
\bibitem{jan}M, Janjic, Determinants and Recurrence Sequences, J. Integer Sequences,
    15 (2012),
Article 12.3.5.
\bibitem{sim}   Yüksel Soykan, Simson Identity of Generalized $m$-step Fibonacci
    Numbers,
Int. J. Adv. Appl. Math. and Mech. 7(2) (2019) 45 – 56.
\end{thebibliography}
\end{document}